\documentclass[11pt]{article}

\usepackage[all]{xy}
\usepackage{amsfonts}
\usepackage{amssymb}
\usepackage{makeidx}
\usepackage{graphicx}
\usepackage{multicol}

\usepackage{amsmath,amssymb,bm}
\usepackage{graphicx}
\usepackage{url}
\usepackage{cite}
\usepackage{enumerate}

\newtheorem{theorem}{Theorem}[section]

\newtheorem{lemma}[theorem]{Lemma}
\newtheorem{proposition}[theorem]{Proposition}

\newtheorem{open problem}[theorem]{Open Problem}

\begin{document}

\title{A Batalin-Vilkovisky Algebra structure on the Hochschild Cohomology of Truncated Polynomials}

\author{Tian Yang }
\maketitle

\begin{abstract}

We calculate the Batalin-Vilkovisky
structure of $HH^*(C^*(\mathbf{K}P^n;R);C^*(\mathbf{K}P^n;R))$ for $
\mathbf{K}=\mathbb{C}$ and $\mathbb{H}$, and $R=\mathbb{Z}$ and any
field; and show that in the special case when
$M=\mathbb{C}P^1=S^2$, and $R=\mathbb{Z}$, this structure can not be
identified with the BV-structure of $\mathbb{H}_*(LS^2;\mathbb{Z})$
computed by Luc Memichi in \cite{menichi2}. However, the induced
Gerstenhaber structures are still identified in this case. Moreover,
according to the work of Y.Felix and J.Thomas
\cite{felix--thomas}, the main result of the present paper
eventually calculates the BV-structure of the rational loop
homology, $\mathbb{H}_*(L\mathbb{C}P^n;\mathbb{Q})$ and
$\mathbb{H}_*(L\mathbb{H}P^n;\mathbb{Q})$, of projective spaces.

\end{abstract}

\section{Introduction}

Let $M$ be a connected, closed oriented manifold of dimension
$d$, and $LM$ the free loop space of $M$. In \cite{chas--sullivan},
\textbf{Theorem 5.4,4.7 and 6.1} Chas and Sullivan defined a
Batalin-Vilkovisky structure on the loop homology $\mathbb{H}_*(LM)$
inducing a Gerstenhaber structure, and a Lie algebra structure on
the string homology $H_*^{S^1}(LM)$. In \cite{cohen--jones},
\textbf{Theorem 3} Cohen and Jones suggested an identification of the
loop homology $\mathbb{H}_*(LM)$ with the Hochschild cohomology
$HH^*(C^*(M);C^*(M))$ as graded algebras. Similar results for
rational and real coefficients were proved by Felix,Thomas and
Vigue-Poirrier in \cite{felix--thomas--vigue-poirrier2} and
 Merkulove in \cite{Mer} respectively. What is interesting is that there is also a way of
defining a Gerstenhaber \cite{gerstenhaber}, and even a BV-structure
on $H^*(C^*(M);C^*(M))$, see \cite{menichi1}, \textbf{Theorem 1.4 a} and \cite{tradler2}, \textbf{Theorem 3.1}, and these two Gerstenhaber,
and BV-structures are expected to be identified.

Also, as a direct consequence of \cite{jones}, \textbf{Theorem A},
there is a nature isomorphism of the string homology of $M$ and the
negative cyclic homology of the chain complex of $M$. In
\cite{menichi1}, \textbf{Theorem 1.4 b}, Luc Menichi defined a Lie
algebra structure on $HC_-^*(C^*(M))$, at least when $M$ is formal;
and it is natural to expect that this Lie algebra structure can be
identified with the string bracket.

For the computational aspect of this theory, there are not so many
results yet. According to the author's knowledge, only for some
really familiar family of manifolds, some partial results of the loop homology
and the string homology are computed, i.e. $S^n$'s, and
$\mathbb{C}P^n$'s. Precisely, for $M=S^n$, in \cite{menichi2}, Luc
Menichi calculates the BV-structures of
$\mathbb{H}_*(LS^2;\mathbb{Z}_2)$ and
$\mathbb{H}_*(LS^2;\mathbb{Z})$.

In \cite{cohen--jones--yan}, Cohen, Jones and Yan developed a
spectral sequence to compute the loop homology of $M=\mathbb{C}P^n$.
However, their paper did not mention the Gerstenhaber and the BV
structures.

For the string homology, in \cite{felix--thomas--vigue-poirrier2},
Felix, Thomas and Vigue-Poirrier has constructed a model for the
string bracket and computed this for $\mathbb{C}P^n$'s with rational
coefficients which turns out to be degenerated.

For the Hochschild cohomology, there are not so many results either.
In \cite{menichi2}, Luc Menichi calculates the BV-structure of
$HH^*(C^*(S^n;\mathbb{Z}_2);C^*(S^n;\mathbb{Z}_2))$, and shows that
this BV-structure can not be identified with the one on
$\mathbb{H}_*(LS^2;\mathbb{Z}_2)$ that he computes, but the induced
Gerstenhaber structures can be identified in this case. In \cite
{westerland2} \textbf{Corollary 4.2}, C.Westerland calculates the
BV-structure of $HH^*(C^*(\mathbb{C}P^n;\mathbb{Z}_2);$
$C^*(\mathbb{C}P^n;\mathbb{Z}_2))$. Both of them are working with the
$\mathbb{Z}_2$ coefficients.

The main result of the present paper describes the Hochschild
cohomology of the cochain complex of $\mathbb{C}P^n$, and also of
$\mathbb{H}P^n$, with coefficients in the integers $\mathbb{Z}$, and
in any field $\mathbf{k}$ with $char(\mathbf{k})=0$ and $\mathbb{F}_p$
with $char(\mathbb{F}_p)=p$, an
arbitrary prime number other than $2$. Precisely,

\begin{quote}
\textbf{Main Theorem:} \textit{Let $A=R[x]/(x^{n+1})$, where
$|x|=2m$ and $R=\mathbb{Z}$ and $\mathbf{k}$; and also
$\mathbb{F}_p$, but in this case we require that $n\neq kp-1$. Then
as a Batalin-Vilkovisky algebra:
$$HH^*(A;A)=R[x,u,t]/(x^{n+1},u^2,(n+1)x^nt,ux^n),$$
where $|x|=-2m,\ |u|=-1\ and\ |t|=2mn+2(m-1)$, and}
\begin{equation*}
\begin{split}
&\Delta(t^kx^l)=0\\
&\Delta(t^kux^l)=(-(k+1)n-k+l)t^kx^l.\\
\end{split}
\end{equation*}
\end{quote}

If $n=kp-1$, we have:
\begin{quote}
\textbf{Theorem \ref {p}} \textit{Let $A=\mathbb{F}_p[x]/(x^{kp})$,
where $|x|=2m$, as a BV-algebra:
$$HH^*(A;A)=\mathbb{F}_p[x,v,t]/(x^{kp},v^2),$$ where $|x|=-2m,\ |v|=2m-1\
and\ |t|=2mn+2(m-1)$, and}
\begin{equation*}
\begin{split}
&\Delta(t^kx^l)=0\\
&\Delta(t^kvx^l)=lt^kx^{l-1}.\\
\end{split}
\end{equation*}
\end{quote}

As an immediate consequence, we have

\textbf{Corollary:} \textit{Let $A=R[x]/(x^{n+1})$, where
$|x|=2m$ and $R=\mathbb{Z}$ and $\mathbf{k}$; and also
$\mathbb{F}_p$, but in this case we require that $n\neq kp-1$. The Gerstenhaber bracket on $HH^*(A;A)$ is determined as follows:}
\begin{equation*}
\begin{split}
&\{t^{k_1}x^{l_1},t^{k_2}x^{l_2}\}=0\\
&\{t^{k_1}x^{l_1},t^{k_2}ux^{l_2}\}=(-k_1n-k_1+l_1)t^{k_1+k_2}x^{l_1+l_2}\\
&\{t^{k_1}ux^{l_1},t^{k_2}ux^{l_2}\}=((k_1+k_2+2)n+k_1+k_2-l_1-l_2)t^{k_1+k_2}ux^{l_1+l_2}.\\
\end{split}
\end{equation*}

\textit{In the case that $A=\mathbb{F}_p[x]/(x^{kp})$,
where $|x|=2m$, the Gerstenhaber bracket is determined as:}
\begin{equation*}
\begin{split}
&\{t^{k_1}x^{l_1},t^{k_2}x^{l_2}\}=0\\
&\{t^{k_1}x^{l_1},t^{k_2}vx^{l_2}\}=l_1t^{k_1+k_2}x^{l_1+l_2-1}\\
&\{t^{k_1}vx^{l_1},t^{k_2}vx^{l_2}\}=-(l_1+l_2)t^{k_1+k_2}vx^{l_1+l_2-1}.\\
\end{split}
\end{equation*}

As a consequence, since $\mathbb{C}P^1=S^2$, there comes:

\textbf{Theorem \ref {s2}} \textit{As Batalin-Vilkovisky algebras,
$\mathbb{H}_*(LS^2;\mathbb{Z})$ is not isomorphic to
$HH^*(C^*(S^2;\mathbb{Z});C^*(S^2;\mathbb{Z}))$, but the induced
Gerstenhaber algebra structures are isomorphic.}

In \cite{felix--thomas},\textbf{Theorem 1} Y.Felix and
J.Thomas proved that the BV-structures on
$\mathbb{H}_*(LM;\mathbf{k})$ and on
$HH^*(C^*(M;\mathbf{k});C^*(M;\mathbf{k}))$ can be naturally
identified. Therefore, the main result of the present paper eventually calculates the BV-structure of the rational loop
homology, $\mathbb{H}_*(L\mathbb{C}P^n;\mathbf{k})$ and
$\mathbb{H}_*(L\mathbb{H}P^n;\mathbf{k})$, of projective spaces.
Precisely, we have:

\textbf{Theorem \ref {lh}} \textit{As BV-algebras:
$$\mathbb{H}_*(L\mathbb{C}P^n;\mathbf{k})=\mathbf{k}[x,u,t]/(x^{n+1},u^2,x^nt,ux^n),$$
with $|x|=-2,\ |u|=-1\ and\ |t|=2n$; and
$$\mathbb{H}_*(L\mathbb{H}P^n;\mathbf{k})=\mathbf{k}[x,u,t]/(x^{n+1},u^2,x^nt,ux^n),$$
with $|x|=-4,\ |u|=-1\ and\ |t|=4n+2$; and in both cases,}
\begin{equation*}
\begin{split}
&\Delta(t^kx^l)=0\\
&\Delta(t^kux^l)=(-(k+1)n-k+l)t^kx^l.\\
\end{split}
\end{equation*}

With the aid of the main theorem, there also comes the calculation
of the negative cyclic cohomology of cochain complex of
$\mathbb{C}P^n$ which was identified by Jones with the string
homology of $\mathbb{C}P^n$ as graded modules. The Lie structure is
also calculated. When the coefficients are in the rationals, the Lie
bracket turns out to be trivial as expected.

\textbf{Theorem \ref {tri}} \textit{As a Lie algebra,}
\begin{equation*}
\begin{split}
HC_-^m(\mathbf{k}[x]/(x^{n+1}))&=\mathbf{k},\ \ m=2q\ \
q=0,1,2,...\\
HC_-^m(\mathbf{k}[x]/(x^{n+1}))&=0,\ \ m=2q+1\ \ q=0,1,2,...\\
HC_-^m(\mathbf{k}[x]/(x^{n+1}))&=0,\ \ m=-2q\ \ q=1,2,...\\
HC_-^m(\mathbf{k}[x]/(x^{n+1}))&=\mathbf{k}^n,\ \ m=-2q+1\ \
q=1,2,...\\
\end{split}
\end{equation*}
\textit{with the Lie bracket $[\_,\_]=0$.}

The method of proving the main theorem is purely homological
algebra. For the graded module structure of $HH^*(A,A)$, it is from
the derived functoriality of the Hochschild cohomology and a
standard $2$-periodic resolution of truncted polynomials

$$(P^*(A),d):\ \ \ 0\rightarrow
A\stackrel{0}{\rightarrow}A\stackrel{(n+1)x^n}{\longrightarrow}A\stackrel{0}{\rightarrow}A\stackrel{(n+1)x^n}{\longrightarrow}\cdots,$$

which was already known by many authors, \cite{loday}
\cite{westerland}.

For the algebra structure, there comes the key step of the argument,
i.e., we carefully constructed a chain map
$$\varphi: P(A)\rightarrow C_*^{bar}(A)$$ from the periodic
resolution to the bar resolution inducing the isomorphism of
homology as graded modules. Therefore, the product can be traced by
the explicit formula of the cup product on the bar complex via
$\varphi^*$. The algebra structure was also known by
\cite{cohen--jones--yan} and the identification in
\cite{cohen--jones}, but the approach in the present paper is more
direct and necessary for getting the BV-structure.

The BV-structure comes from the checking of a relatively manageable
formula for the $\Delta$-operator. This formula for $\Delta$ makes
use of the Poincare dual basis. The Gerstenhaber structure is a
direct consequence of the BV-structure.

The paper is organized as follows. In section $2$, the Hochschild
cohomology of differential graded algebras is reviewed. The explicit
formula for the cup product and the $\Delta$-operator is obtained.
Section $3$ concentrates on the proof of the main theorem. Section
$4$ consists of some consequences of the main theorem.

\section{An explicit formula of the $\Delta$-operator}
In this section an explicit formula for the $\Delta$-operator which
will play a key role in the later calculation is obtained.

According to \cite{menichi1}, \textbf{Theorem 1.6} when $A$ has
certain symmetry, such as the Poincar$\acute{e}$
duality, the Gerstenhaber structure can be extended to a
BV-structure as in the non-$DG$ case, at least for the special case
when the differential $d$ is trivial. To this end, we should employ the
Connes' boundary-operator with a slight modification of the sign.

\textbf{Connes' boundary-operator for $DG$-algebra:} Let $A$ be a
$DG$-algebra, the Connes' boundary operator
$B:A\otimes(s\overline{A})^{\otimes n}\rightarrow
A\otimes(s\overline{A})^{\otimes n+1}$ on the Hochschild chain
complex of $A$ with coefficients in itself, $C_*(A;A)$, is defined
by
$$B(a_0\otimes(a_1,...,a_n))=\sum_{i=0}^n(-1)^{\sum_{k=0}^{i-1}|sa_k|\sum_{k=i}^n|sa_k|}1\otimes(a_i,...,a_n,a_0,a_1,...,a_{i-1}).$$

To get a BV-structure on the Hochschild cohomology $HH^*(A)$, we
should use a sequence of isomorphisms coming from the duality of the
algebra $A$ and the adjunction between the tensor and $Hom$
functors. The duality gives an isomorphism
$Hom(T(s\overline{A}),A)\cong Hom(T(s\overline{A}),A^*)$ as chain
complexes, and the adjunction gives $Hom(T(s\overline{A}),A^*)\cong
Hom(A\otimes T(s\overline{A}),k)$, where the right hand side is the dual of the Hochschild chain complex. Therefore, we can put the
dual of the Connes' boundary-operator onto
$Hom(T(s\overline{A}),A)$, the Hochschild cochain complex, via this
string of isomorphisms. Chasing each of the isomorphisms carefully,
we have:

\begin{proposition}\label{Delta}
The operator $\Delta: Hom(s\overline{A}^{\otimes n+1},A)\rightarrow
Hom(s\overline{A}^{\otimes n},A)$ is given by
\begin{equation*}
\begin{split}
&\Delta(f)((a_1,...,a_n))\\
&=\sum_{j=1}^N(-1)^{|f|+|a^j|\sum_{k=1}^n|sa_k|}\sum_{i=0}^n(-1)^{(|sa_0|+\sum_{k=1}^{i-1}|sa_k|)\sum_{k=i}^n|sa_k|}
<1,f((a_i,...,a_n,a_0,a_1,...,a_{i-1}))>{a^j}^*,\\
\end{split}
\end{equation*}

for all $f\in Hom(s\overline{A}^{\otimes n+1},A)$.
\end{proposition}

where $a_0=a^j$ in this formula. That is the explicit formula of
$\Delta$-operator with which we can do some concrete calculation;
and by Luc Menichi's theorem \cite{menichi1},\textbf{Theorem 1.6},
this $\Delta$-operator induces a BV structure on $HH^*(A;A)$ which
induces the Gerstenhaber structure of
$HH^*(A;A)$.

\section{The main theorem}
In this section, we prove the main theorem of the present paper.
For convenience, $A$ will always mean
$\mathbb{Z}[x]/(x^{n+1})$ throughout this section.

\textbf{Main Theorem:} \textit{Let $A=\mathbb{Z}[x]/(x^{n+1})$,
where $|x|=2$, then as a Batalin-Vilkovisky algebra:
$$HH^*(A;A)=\mathbb{Z}[x,u,t]/(x^{n+1},u^2,(n+1)x^nt,ux^n),$$
where $|x|=-2,\ |u|=-1\ and\ |t|=2n$, and}
\begin{equation*}
\begin{split}
&\Delta(t^kx^l)=0\\
&\Delta(t^kux^l)=(-(k+1)n-k+l)t^kx^l.\\
\end{split}
\end{equation*}

\subsection{The graded $R$-module structure} 

By the $Tor$ interpretation of $HH^*(A;A)$, it is much more convenient to
take a simple projective resolution other than the bar resolution.
We have:

\begin{proposition}The following
$$P_*(A):\ \ \ \ \cdots\rightarrow \Sigma^{2(n+1)}(A\otimes
A)\stackrel{y^n+y^{n-1}z+\cdots
+z^n}{\longrightarrow}\Sigma^2(A\otimes
A)\stackrel{y-z}{\rightarrow}A\otimes
A\stackrel{\mu}{\rightarrow}A\rightarrow 0$$ gives a 2-periodical
resolution of $A$ as $A\otimes A$-module with
$P_{2k}(A)=\Sigma^{2k(n+1)} A\otimes A$ with $d_{2k}(a\otimes
b)=y^na\otimes b+y^{n-1}a\otimes zb+\cdots+a\otimes z^nb$ and
$P_{2k+1}=\Sigma^{2k(n+1)+2} A\otimes A$ with $d_{2k+1}(a\otimes
b)=ya\otimes b-a\otimes zb$, where $\Sigma$ is the degree increasing
operator which makes the differential of degree $0$ ,and $y,z$ are
generators of the copies of $A$ in $A\otimes A$ respectively.
\end{proposition}

\textit{Proof}: See \cite{loday} and \cite{weibel}. \hfill$\square$

After taking $Hom_{A\otimes A}(\_,A)$ of this periodical complex
$P_*(A)$, we get the following periodical cochain complex
$$(P^*(A),d):\ \ \ 0\rightarrow
A\stackrel{0}{\rightarrow}A\stackrel{(n+1)x^n}{\longrightarrow}A\stackrel{0}{\rightarrow}A\stackrel{(n+1)x^n}{\longrightarrow}\cdots$$
by chasing the differential. Therefore, a direct computation tells us:

\begin{proposition} As abelian groups,
\begin{equation*}
\begin{split}
&HH^0(A;A)=A=\mathbb{Z}[x]/(x^{n+1})\\
&HH^{2q-1}(A;A)=ker(n+1)x^n=\overline{A}\cong
\mathbb{Z}^n\\
&HH^{2q}(A;A)=A/((n+1)x^n)\cong \mathbb{Z}^n\oplus
\mathbb{Z}_{n+1}.\\
\end{split}
\end{equation*}
\end{proposition}

Where $\mathbb{Z}_m$ stands for $\mathbb{Z}/m\mathbb{Z}$ for typing
convenience.

\subsection{Key step: the construction of $\varphi$}

To get the algebra structure of $HH^*(A;A)$, we go back to
$Hom(Ts\overline{A},A)$ on which the cup product is defined. To do this, we define $A\otimes A$-module maps
$$\varphi_*: P_*(A)\rightarrow \overline{C}_*^{bar}(A),$$ 
with
$$\varphi_{2q}:\Sigma^{2q(n+1)}A\otimes A\rightarrow A\otimes (s\overline{A})^{\otimes 2q}\otimes
A$$ defined by
\begin{equation}\label{phi1}
\varphi_{2q}(1\otimes1)=\sum1[x^{n-a_1}|x|
x^{n-a_2}|x|\cdots|x^{n-a_q}|x]x^{\sum_{k=1}^qa_k},
\end{equation}
and

$$\varphi_{2q+1}:\Sigma^{2q(n+1)+2}A\otimes A\rightarrow
A\otimes(s\overline{A})^{\otimes 2q+1}\otimes A$$ defined by
\begin{equation}\label{phi2}
\varphi_{2q+1}(1\otimes1)=\sum1[x| x^{n-a_1}| x|
x^{n-a_2}|x|\cdots|x^{n-a_q}|x]x^{\sum_{k=1}^qa_k},
\end{equation}
where the sum is taken over $0\leqslant a_k<n$, $k=1,2,...,q$ and
$\sum_{k=1}^qa_k\leqslant n$. We have the following key lemma:

\begin{lemma} $\varphi^*=Hom_{A\otimes A}(\varphi,A):C^*(A;A)\rightarrow
P^*(A)$ is a cochain map.
\end{lemma}

\textit{Proof:} Since
$\Phi^*:(Hom_k(T(s\overline{A}),A),\beta)\rightarrow
(Hom_{A^e}(A\otimes T(s\overline{A})\otimes A,A),\beta')$ is an
isomophism of cochain complex, it suffices to show $d\circ
(\varphi^*\circ\Phi^*)=(\varphi^*\circ\Phi^*)\circ\beta:Hom_k(T(s\overline{A}),A)\rightarrow
P^*(A)$, i.e. for any $f\in Hom_k(s\overline{A}^{\otimes m},A)$,
$$(d\circ
\varphi_m^*\circ\Phi^*)(f)(1\otimes1)=(\varphi_{m+1}^*\circ\Phi^*\circ\beta)(f)(1\otimes1),$$
where $1\otimes1\in P_m(A)$.

For $m=2q$, we have on one side $(d_{2q}\circ \varphi_{2q}^*)(f)=0$,
since $d_{2q}=0$. The calculation for the other side is more complicated. We have

\begin{equation*}
\begin{split}
\varphi_{2q+1}^*\circ\Phi^*(\beta(f))(1\otimes1)&=\Phi^*(\beta(f))(\phi_{2q+1}(1\otimes
1))\\
&=\beta'(\Phi^*(f))(\sum_{0\leqslant a_k<n}1[x| x^{n-a_1}| x|\cdots|
x^{n-a_{q}}| x] x^{\sum_{k=1}^{q}a_k})\\
&=\sum_{0\leqslant a_k<n}x^{\sum_{k=1}^{q}a_k}\beta(f)(x, x^{n-a_1},
x,\cdots, x^{n-a_{q}},x),\\
\end{split}
\end{equation*}

the last equality is from the commutativity of $A$ and that the
degree of $x$ is even. Now as a typical term,

\begin{equation*}
\begin{split}
T&=x^{\sum_{k=1}^{q}a_k}\beta(f)(x, x^{n-a_1}, x,\cdots,
x^{n-a_{q}}, x)\\
&=x^{\sum_{k=1}^{q}a_k+1}f(x^{n-a_1}, x,\cdots, x^{n-a_{q}}, x)\\
&-\sum_{k=1}^qx^{\sum_{k=1}^{q}a_k}f(x, x^{n-a_1},\cdots,
x^{n-a_{i-1}}, x^{n-a_i+1}, x,\cdots, x^{n-a_q},x)\\
&+\sum_{k=1}^qx^{\sum_{k=1}^{q}a_k}f(x, x^{n-a_1},\cdots\, x,
x^{n-a_i+1}, x^{n-a_{i+1}},\cdots, x^{n-a_q}, x)\\
&-x^{\sum_{k=1}^{q}a_k+1}f(x, x^{n-a_1},\cdots, x, x^{n-a_q}).\\
\end{split}
\end{equation*}

Note that the other term $$x^{\sum_{k=1}^{q}a_k+1}\beta(f)(x,
x^{n-a_1-1}, x,\cdots, x^{n-a_q}, x)$$ has a term 
$$-x^{\sum_{k=1}^{q}a_k+1}f(x^{n-a_1}, x,\cdots, x^{n-a_{q}},
x)$$ as part of it, which cancels the first term
$$x^{\sum_{k=1}^{q}a_k+1}f(x^{n-a_1}, x,\cdots,
x^{n-a_{q}}, x)$$ in
$$T=x^{\sum_{k=1}^{q}a_k}\beta(f)(x, x^{n-a_1},
x,\cdots, x^{n-a_{q}}, x);$$ and

$$x^{\sum_{k=1}^{q}a_k+1}\beta(f)(x,
x^{n-a_1}, x,\cdots, x^{n-a_{q}-1}, x)$$ has term
$$x^{\sum_{k=1}^{q}a_k+1}f(x, x^{n-a_1},\cdots,
x, x^{n-a_q})$$ which cancels the last term
$$-x^{\sum_{k=1}^{q}a_k+1}f(x, x^{n-a_1},\cdots,
x, x^{n-a_q})$$ in $$T=x^{\sum_{k=1}^{q}a_k}\beta(f)(x, x^{n-a_1},
x,\cdots, x^{n-a_{q}}, x);$$

also,
$$x^{\sum_{k=1}^{q}a_k}\beta(f)(x, x^{n-a_1}, x,\cdots, x^{n-a_i+1}, x, x^{n-a_{i+1}-1},\cdots, x^{n-a_{q}}, x)$$
have terms $$-x^{\sum_{k=1}^{q}a_k}f(x, x^{n-a_1},\cdots, x,
x^{n-a_i+1}, x^{n-a_{i+1}},\cdots, x^{n-a_q}, x)$$ which cancel
terms
$$x^{\sum_{k=1}^{q}a_k}f(x, x^{n-a_1},\cdots,
x, x^{n-a_i+1}, x^{n-a_{i+1}},\cdots, x^{n-a_q}, x)$$ in
$$T=x^{\sum_{k=1}^{q}a_k}\beta(f)(x, x^{n-a_1},
x,\cdots, x^{n-a_{q}}, x);$$ and
$$x^{\sum_{k=1}^{q}a_k}\beta(f)(x, x^{n-a_1}, x,\cdots, x^{n-a_{i-1}-1}, x, x^{n-a_i+1},\cdots, x^{n-a_q}, x)$$
give rise to terms $$x^{\sum_{k=1}^{q}a_k}f(x, x^{n-a_1},\cdots,
x^{n-a_{i-1}}, x^{n-a_i+1}, x,\cdots, x^{n-a_q}, x)$$ which cancel
$$-x^{\sum_{k=1}^{q}a_k}f(x, x^{n-a_1},\cdots,
x^{n-a_{i-1}}, x^{n-a_i+1}, x,\cdots, x^{n-a_q},x)$$ in
$$T=x^{\sum_{k=1}^{q}a_k}\beta(f)(x, x^{n-a_1},
x,\cdots, x^{n-a_{q}}, x).$$ Therefore, for each $a_k$,
$0<a_k<n-1$,
$$x^{\sum_{k=1}^{q}a_k}\beta(f)(x, x^{n-a_1},
x,\cdots, x^{n-a_{q}}, x)$$ is canceled by other terms in
$\varphi_{2q+1}^*(\beta(f))(1\otimes1).$

Now, we only need to focus on the edge effects, i.e. those terms
with $a_k=n-1$ or $0$, and those not get canceled in the way
above. When $a_1=n-1$,
$$x^{n-1}\beta(f)(x,x, x,x^n,\cdots,
x^n,x)$$ gives rise to one term $$x^nf(x, x, x^n,\cdots, x^n, x)$$
which is unable to be canceled by the way above, but it gets
canceled by a term that also could not be canceled in the way above
coming from
$$x^n\beta(f)(x, x, x, x^{n-1},\cdots,
x^n, x).$$

When $a_k=n-1,\ 0<k<q,$ the terms $$-x^nf(x, x^n, x,\cdots, x, x,
x^n,\cdots, x^n, x)$$ coming from
$$x^n\beta(f)(x, x^n, x,\cdots, x,
x, x, x^{n-1},\cdots, x^n, x)$$ cancels those corresponding positive
ones given by
$$x^n\beta(f)(x, x^n, x,\cdots,
x^{n-1}, x, x, x,\cdots,x^n, x).$$

When $a_q=n-1$, $$-x^nf(x,x^n, x,\cdots, x^n, x, x)$$ given by
$$x^{n-1}\beta(f)(x, x^n,\cdots, x, x,
x)$$ cancels the corresponding positive one given by
$$x^n\beta(f)(x, x^n, x,\cdots,
x^{n-1}, x, x, x).$$

Finally, when $a_k=0,\ 0\leqslant k \leqslant q$, the terms which
could not get canceled vanish automatically by the degree reason.
Therefore, all the terms are mutually canceled, and
$\varphi_{2q+1}^*(\beta(f))(1\otimes1)$ vanishes as expected.

Now for $m=2q-1$, we have on one side

\begin{equation*}
\begin{split}
(d_{2q-1}\circ
\varphi_{2q-1}^*\circ\Phi^*)(f)(1\otimes1)&=(n+1)x^n\Phi^*(f)(\varphi(1\otimes1))\\
&=(n+1)x^n\hat{f}(\sum_{0\leqslant a_k<n}1[ x| x^{n-a_1}| x| \cdots|
x^{n-a_{q-1}}|x] x^{\sum_{k=1}^{q-1}a_k})\\
&=(n+1)x^n\sum_{0\leqslant a_k<n}f(x, x^{n-a_1}, x,\cdots,
x^{n-a_{q-1}}, x)x^{\sum_{k=1}^{q-1}a_k}\\
&=(n+1)x^nf(x, x^n, x,\cdots, x^n, x).\\
\end{split}
\end{equation*}

For the other side,  we have
\begin{equation*}
\begin{split}
\varphi_{2q}^*\circ\Phi^*(\beta(f))(1\otimes1)&=\beta'(\Phi^*(f))(\varphi_{2q}^*(1\otimes1))\\
&=\sum_{0\leqslant a_k<n}x^{\sum_{k=1}^q a_k}\beta(f)(x^{n-a_1},
x,\cdots, x^{n-a_q}, x),\\
\end{split}
\end{equation*}

and
\begin{equation*}
\begin{split}
&x^{\sum_{k=1}^q a_k}\beta(f)(x^{n-a_1}, x,\cdots, x^{n-a_q},
x)\\
=&x^{n+\sum_{i=2}^qa_i}f(x, x^{n-a_2},\cdots, x^{n-a_q}, x)\\
-&\sum_{i=1}^qx^{\sum_{k=1}^qa_k}f(x^{n-a_1}, x,\cdots, x^{n-a_i+1},
x^{n-a_{i+1}}, x,\cdots, x^{n-a_q}, x)\\
+&\sum_{i=1}^{q-1}x^{\sum_{k=1}^qa_k}f(x^{n-a_1}, x,\cdots,
x^{n-a_i}, x^{n-a_{i+1}+1}, x,\cdots, x^{n-a_q}, x)\\
+&x^{\sum_{k=1}^q+1}f(x^{n-a_1}, x,\cdots, x, x^{n-a_q}).\\
\end{split}
\end{equation*}

Similar to the previous case, most of the terms cancel each other,
and the only exceptions are: $$x^n\beta(f)(x, x, x^{n-1},\cdots,
x^n, x)$$ gives rise to one term
$$x^nf(x, x^n, x,\cdots, x^n, x)$$
that cannot be canceled; and each $$x^i\beta(f)(x^{n-i}, x,\cdots,
x^n, x),\ i=0,1,...,n-1,$$ gives rise to the term
$$x^nf(x,x^n, x,\cdots, x^n, x)$$
that survives from canceling. Therefore,
$$\varphi_{2q}^*(\beta(f))(1\otimes1)=(n+1)x^nf(x, x^n, x,\cdots,
x^n,x)$$ as expected. This completes the proof of the lemma.
\hfill$\square$

\subsection{The graded commutative algebra structure} Note that,
there are three distinguished elements $\bar{x},\bar{u}$ and
$\bar{t}$ in $C^*(A;A)$, where
\begin{equation*}
\begin{split}
&\bar{x}\in Hom_k(k,A)\ \ with\ \ \bar{x}(1)=x,\\
&\bar{u}\in Hom_k(s\overline{A},A)\ \ with\ \ \bar{u}(x^i)=ix^i,\ \
\text{and}\\
&\bar{t}\in Hom_k(s\overline{A}\otimes s\overline{A},A)\ \ with\ \
\bar{t}(x^i, x^j)=x^{i+j-(n+1)},\\
\end{split}
\end{equation*}

with $|\bar{x}|=-2$, $|\bar{u}|=-1$ and $|\bar{t}|=2n$.

The significance of these elements is that they represent
non-trivial cohomology classes, and moreover, give rise to the
generators of $HH^*(A;A)$ as a graded commutative algebra.
\begin{lemma} The elements $\bar{x},\bar{u}$
and $\bar{t}$ in $C^*(A;A)$ are cocycles but not coboundaries, hence
present non-trivial cohomology classes in $HH^*(A;A)$.
\end{lemma}

\textit{Proof:} Since
\begin{equation*}
\begin{split}
\beta(\bar{x})(x^i)&=x^i\bar{x}(1)-\bar{x}(1)x^i\\
&=x^ix-xx^i\\
&=0,\ 0<i\leqslant n;\\
\end{split}
\end{equation*}

and the cochain complex is non-negative, $\bar{x}$ represents a
non-trivial class of $HH^0(A,A)$.

Also, if $i+j\leqslant n,$
\begin{equation*}
\begin{split}
\beta(\bar{u})(x^i\otimes
x^j)&=x^i\bar{u}(x^j)-\bar{u}(x^{i+j})+\bar{u}(x^i)x^j\\
&=jx^{i+j}-(i+j)x^{i+j}+ix^{i+j}\\
&=0;\\
\end{split}
\end{equation*}

and if $i+j>n,$
$$\beta(\bar{u})(x^i\otimes
x^j)=x^i\bar{u}(x^j)+\bar{u}(x^i)x^j=0,\ \ since\ x^{i+j}=0.$$ Hence
$\bar{u}$ is a cocycle. $\bar{u}$ is not a coboundary, since for any
$g\in Hom_k(k,A)$,
$$\beta(g)(x^i)=x^ig(1)-g(1)x^i=0.$$ Therefore, $\bar{u}$ represents
a non-trivial class in $HH^1(A;A)$.

Lastly, for $\bar{t}$, we have $$\beta(\bar{t})(x^k, x^i,
x^j)=x^k\bar{t}(x^i, x^j)-\bar{t}(x^{k+i}, x^j)+\bar{t}(x^k,
x^{i+j})-\bar{t}(x^k, x^i)x^j.$$ At this point, we do a case by case discussion. When
$k+i<n+1$, $i+j<n+1$,
\begin{equation*}
\begin{split}
\beta(\bar{t})(x^k, x^i, x^j)&=-\bar{t}(x^{k+i}, x^j)+\bar{t}(x^k,
x^{i+j})\\
&=-x^{k+i+j-(n+1)}+x^{k+i+j-(n+1)}=0;\\
\end{split}
\end{equation*}

when $k+i<n+1$, $i+j\geqslant n+1$ (similarly, $k+i\geqslant n+1$,
$i+j<n+1$), $x^{i+j}=0$ so
\begin{equation*}
\begin{split}
\beta(\bar{t})(x^k,x^i,x^j)&=x^k\bar{t}(x^i, x^j)-\bar{t}(x^{k+i},
x^j)\\
&=x^{k+i+j-(n+1)}-x^{k+i+j-(n+1)}=0;\\
\end{split}
\end{equation*}

and when $k+i\geqslant n+1$, $i+j\geqslant n+1$,
\begin{equation*}
\begin{split}
\beta(\bar{t})(x^k, x^i, x^j)&=x^k\bar{t}(x^i, x^j)-\bar{t}(x^k,
x^i)x^j\\
&=x^{k+i+j-(n+1)}-x^{k+i+j-(n+1)}=0.\\
\end{split}
\end{equation*}

Therefore, $\bar{t}$ is a cocycle.

To see that it is not a coboundary, we have, $\forall g\in
Hom_k(s\overline{A},A)$, $$\beta(g)(x^n,
x^n)=x^ng(x^n)+g(x^n)x^n=2x^ng(x^n)$$ which cannot be equal to
$\bar{t}(x^n, x^n)=x^{n-1}$ for dimensional reasons. Therefore,
$\bar{t}$ is not a coboundary, hence represents a non-trivial
cohomology class in $HH^2(A;A)$.\hfill$\square$
\begin{lemma}

\begin{enumerate}
\item As a commutative graded algebra, $HH^*(A;A)$ is generated by
$\bar{x}$, $\bar{u}$ and $\bar{t}$.

\item $\varphi^*$ induces an isomorphism on homology.
\end{enumerate}
\end{lemma}

\textit{Proof:} By (\ref{phi1}) and (\ref{phi2}),
$\varphi^*(\bar{x})(1\otimes 1)=\bar{x}(1)=x,$ so
$\varphi^*(\bar{x})=x\in A\cong HH^0(A;A);$
$\varphi^*(\bar{u})(1\otimes 1)=\bar{u}(x)=x,$ so
$\varphi^*(\bar{u})=x\in ker(n+1)x^n\cong HH^1(A;A);$ and
$\varphi^*(\bar{t})(1\otimes1)=\sum_{k=0}^{n-1}\bar{t}(x^{n-k},
x)x^k=1,$ hence $\varphi^*(\bar{t})=1\in A/((n+1)x^n)\cong
HH^2(A;A).$ Moreover, for $\bar{t}^q\bar{x}^l\in
Hom_k(s\overline{A}^{2q},A)$, we have
\begin{equation*}
\begin{split}
\varphi^*(\bar{t}^q\bar{x}^l)(1\otimes 1)&=\sum_{0\leqslant a_k<
n}\bar{t}^q\bar{x}^l(x^{n-a_1}, x, x^{n-a_2}, x,\cdots, x^{n-a_q},
x)x^{\sum_{k=1}^qa_k}\\
&=\bar{t}^q\bar{x}^l(x^n, x,\cdots, x^n, x)\\
&=x^l,\\
\end{split}
\end{equation*}

hence
$$\varphi^*(\bar{t}^q\bar{x}^l)=x^l\in A/((n+1)x^n)\cong
HH^{2q}(A;A);$$ and for $\bar{t}^q\bar{u}\bar{x}^l\in
Hom_k(s\overline{A}^{2q+1},A)$, we have
\begin{equation*}
\begin{split}
\varphi^*(\bar{t}^q\bar{u}\bar{x}^l)(1\otimes
1)&=\sum_{0\leqslant a_k< n}\bar{t}^q\bar{u}\bar{x}^l(x, x^{n-a_1},
x, x^{n-a_2}, x,\cdots, x^{n-a_q},
x)x^{\sum_{k=1}^qa_k}\\
&=\bar{t}^q\bar{u}\bar{x}^l(x, x^n, x,\cdots, x^n,
x)\\
&=x^{l+1},\\
\end{split}
\end{equation*}

hence
$$\varphi^*(\bar{t}^q\bar{u}\bar{x}^l)=x^{l+1}\in ker(n+1)x^n\cong
HH^{2q+1}(A;A).$$

Therefore, $\varphi^*$ induces a surjection on homology, Since
we are computing a value of $Tor$
which is finitely generated using two projective resolutions, $\varphi^*$ induces an
isomorphism on homology. Therefore, $\bar{x}$, $\bar{u}$ and
$\bar{t}$ generate $HH^*(A;A)$ as an algebra.\hfill$\square$

\begin{proposition} As a commutative graded
algebra,$$HH^*(A;A)=\mathbb{Z}[x,u,t]/(x^{n+1},u^2,(n+1)x^nt,ux^n).$$
\end{proposition}

\textit{Proof:} Let $x$, $u$ and $t$ denote the cohomology classes of $\bar{x}$, $\bar{u}$ and $\bar{t}$
respectively. The relation between them is straightforward on
$P^*(A)$.

Precisely, consider $x=\varphi^*(\bar{x})\in P^0(A)=Hom_{A\otimes
A}(A\otimes A,A)$. We have
$x^{n+1}(1\otimes1)=(x(1\otimes1))^{n+1}=x^{n+1}=0\in A$, so
$x^{n+1}=0\in P^0(A).$ Therefore, $\bar{x}^{n+1}=0\in HH^*(A;A).$
The class $\bar{u}\bar{x}^n=0$ for the same reason. $\bar{u}^2=0,$
since $u^2(1\otimes
1)=\varphi^2(\bar{u}^2)(1\otimes1)=\bar{u}^2(\varphi_2(1\otimes1))=\sum_{k=0}^{n-1}\bar{u}^2(x^{n-k},x)x^k=0\in
A$. The element $(n+1)x^nt=\varphi^*((n+1)\bar{x}^n\bar{t})\in
P^2(A)$ is a coboundary, i.e. the image of $1\in P^1(A)$, hence
vanishes in the homology $HH^*(A;A).$ All the other relations are
generated by those four above.\hfill$\square$

The algebra structure was also known to \cite{cohen--jones--yan} and
the identification in \cite{cohen--jones}; but the approach in the
present paper is more direct and necessary for getting the
BV-structure.

\subsection{The BV-structure}

We are now left to find the BV-structure. Since in a BV-algebra,
we have equation \cite{getzler}:
\begin{equation*}
\begin{split}
\Delta(abc)=&\Delta(ab)c+(-1)^{|a|}\Delta(bc)+(-1)^{(|a|-1)|b|}b\Delta(ac)\\
&-\Delta(a)bc-(-1)^{|a|}a(\Delta(b))c-(-1)^{|a|+|b|}ab\Delta(c),\\
\end{split}
\end{equation*}

it suffice to find the value of
$\Delta(x)$,$\Delta(x^2)$,$\Delta(u)$,$\Delta(t)$,$\Delta(t^2)$,$\Delta(tx)$,$\Delta(tu)$
and $\Delta(ux)$; and these values determine $\Delta$ via the
equation. This leads us to the following calculation:

$\Delta(x^k)=0$ by the degree reason. By proposition \ref{Delta}, we
take $\{x^k\}_{0\leqslant k\leqslant n}$ as basis (hence
$\{x^{n-k}\}_{0\leqslant k\leqslant n}$ as dual basis), and have
$$\Delta(u)(1)=\sum_{k=0}^n-<1,u(x^k)>x^{n-k}=-<1,nx^n>1=-n,$$
hence $\Delta(u)=-n.$
$$\Delta(ux)(1)=\sum_{k=0}^n-<1,ux(x^k)>x^{n-k}=-<1,(n-1)x^n>x=-(n-1)x,$$
hence $\Delta(ux)=-(n-1)x.$
$$\Delta(t)(x^i)=\sum_{k=0}^n<1,t(x^k,
x^i)>x^{n-k}-\sum_{k=0}^n<1,t(x^i, x^k)>x^{n-k}=0,$$ since
$t(x^i,x^k)=t(x^k, x^i)$; hence $\Delta(t)=0$.
\begin{equation*}
\begin{split}
\Delta(t^2)(x^i, x^j, x^h)&=\sum_{k=0}^n<1,t^2(x^k, x^i, x^j,
x^h)>x^{n-k}\\
&-\sum_{k=0}^n<1,t^2(x^i, x^j, x^h,
x^k)>x^{n-k}\\
&+\sum_{k=0}^n<1,t^2(x^j, x^h, x^k,
x^i)>x^{n-k}\\
&-\sum_{k=0}^n<1,t^2(x^h, x^k, x^i, x^j)>x^{n-k}.\\
\end{split}
\end{equation*}

To make $<1,\_>$ non-zero, $i+j+k+h-2(n+1)$ should be equal to $n$,
i.e. $i+j+k+h$=3n+2. Therefore, all of $\{k+i,i+j,j+h,h+k\}$ should
be $\geqslant n+1,$ so
\begin{equation*}
\begin{split}
\Delta(t^2)(x^i, x^j,
x^h)&=x^{i+j+h-2(n+1)}-x^{i+j+h-2(n+1)}+x^{i+j+h-2(n+1)}-x^{i+j+h-2(n+1)}\\
&=0,\\
\end{split}
\end{equation*}
hence $\Delta(t^2)=0.$
\begin{equation*}
\begin{split}
\Delta(tx)(x^i)&=\sum_{k=0}^n<1,tx(x^k,
x^i)>x^{n-k}-\sum_{k=0}^n<1,tx(x^i, x^k)>x^{n-k}\\
&=0,\\
\end{split}
\end{equation*}

hence $\Delta(tx)=0.$
$$\Delta(tu)(x^i, x^j)=-\sum_{k=0}^n<1,tu(x^k,
x^i, x^j)>x^{n-k}$$$$-\sum_{k=0}^n<1,tu(x^i, x^j,
x^k)>x^{n-k}-\sum_{k=0}^n<1,tu(x^j, x^k, x^i)>x^{n-k}.$$ When
$i+j<n+1$, $\forall k$, $<1,\_>$ will vanish; hence
$\Delta(tu)(x^i,x^j)=0$ in this case. When $i+j\geqslant n+1$, to
make $<1,\_>$ non-vanished, $k$ should be equal to $k_0=2n+1-(i+j)$,
then $k_0+i=2n+1-(i+j)+i=2n+1-j\geqslant n+1$, and similarly,
$j+k_0\geqslant n+1$. Therefore, in this case,
\begin{equation*}
\begin{split}
&\Delta(tu)(x^i, x^j)\\
&=-<1,tu(x^{k_0}, x^i,
x^j)>x^{n-k_0}-<1,tu(x^i, x^j, x^{k_0})>x^{n-k_0}-<1,tu(x^j,
x^{k_0},
x^i)>x^{n-k_0}\\
&=-jx^{(i+j)-(n+1)}-(2n-(i+j))x^{(i+j)-(n+1)}-ix^{(i+j)-(n+1)}\\
&=-(2n+1)x^{(i+j)-(n+1)}.\\
\end{split}
\end{equation*}

Therefore, $\Delta(tu)=-(2n+1)t.$

Now we have all the data we need to determine $\Delta$. An induction on the powers of $x$ and $t$ tells us that:
\begin{equation*}
\begin{split}
&\Delta(t^kx^l)=0\\
&\Delta(t^kux^l)=(-(k+1)n-k+l)t^kx^l.\\
\end{split}
\end{equation*}

That completes the proof of the main theorem.\hfill$\square$

Note that there is no essential difference between $|x|=2$ and
$|x|=2m$. In the later case, we go through the same argument word by word, and have

\begin{theorem}\label{H}

Let $A=\mathbb{Z}[x]/(x^{n+1})$, where $|x|=2m$, then as a
Batalin-Vilkovisky algebra:
$$HH^*(A;A)=\mathbb{Z}[x,u,t]/(x^{n+1},u^2,(n+1)x^nt,ux^n),$$
where $|x|=-2m,\ |u|=-1\ and\ |t|=2mn+2(m-1)$, and
\begin{equation*}
\begin{split}
&\Delta(t^kx^l)=0\\
&\Delta(t^kux^l)=(-(k+1)n-k+l)t^kx^l.\\
\end{split}
\end{equation*}
\end{theorem}

\section{Consequences of the main theorem}

In this section, we give some consequences of the main theorem.

\subsection{Consequence 1}

\begin{theorem}\label{s2}
As Batalin-Vilkovisky algebras, $\mathbb{H}_*(LS^2;\mathbb{Z})
\ncong HH^*(C^*(S^2;\mathbb{Z});(C^*(S^2;\mathbb{Z}))$, but the
induced Gerstenhaber algebra structures are isomorphic.
\end{theorem}

\textit{Proof:} In the special case when $n=1$,
$A=\mathbb{Z}[x]/(x^2)$, the main theorem gives: as a BV-algebra,
$$HH^*(A)=\mathbb{Z}[x,u,t]/(x^2,u^2,2xt,ux),$$
$$\Delta(t^kx)=0,\ \ \Delta(t^ku)=-(2k+1)t^k.$$

Recall Luc Menichi's result in \cite{menichi2}, \textbf{Theorem 25},
as a BV-algebra,
$$\mathbb{H}_*(LS^2;\mathbb{Z})=\mathbb{Z}[a,b,v]/(a^2,b^2,ab,2av),$$
$$\Delta(v^ka)=0,\ \ \Delta(v^kb)=(2k+1)v^k+av^{k+1}.$$

This is not isomorphic to
$HH^*(\mathbb{Z}[x]/(x^2);\mathbb{Z}[x]/(x^2))$ as BV-algebras,
since by dimensional reason, any isomorphisms of algebras
 $$\Phi: HH^*(\mathbb{Z}[x]/(x^2);\mathbb{Z}[x]/(x^2))\rightarrow
 \mathbb{H}_*(LS^2)$$
must maps $x$ to $\pm a$, $u$ to $\pm b$ and $t$ to $\pm v$ or $\pm
v+av^2$, but none of them could be a BV-isomorphism. Since if
$\Phi(t)=\pm v$, then
$$\Phi\circ\Delta(ut^k)=\pm(2k+1)v^k,$$ but
$$\Delta\circ\Phi(ut^k)=\pm(2k+1)v^k+av^{k+1}\neq\Phi\circ\Delta(ut^k).$$
If $\Phi(t)=\pm v+av^2$, then by induction
$$\Phi\circ\Delta(ut^k)=\Phi(-(2k+1)t^k)=\pm(2k+1)v^k+kav^{k+1},$$
but $$\Delta\circ\Phi(ut^k)=\Delta(\pm
bv^k)=\pm(2k+1)v^k+av^{k+1}\neq\Phi\circ\Delta(ut^k)$$ when $k$ is
even. So there is no possibility for $\Phi$ to be a BV-isomorphism.
However, if we let
$$\Phi:HH^*(\mathbb{Z}[x]/(x^2);\mathbb{Z}[x]/(x^2))\rightarrow\mathbb{H}_*(LS^2)$$
given by $$\Phi(x)=a,\ \
 \Phi(u)=-b\ \ and\ \Phi(t)=v;$$ checking the formula
 $$\{x,y\}=(-1)^{|x|}\Delta(xy)-(-1)^{|x|}(\Delta x)y-a(\Delta b)$$ shows that $\Phi$ preservers the
 Gerstenhaber structure. \hfill$\square$

\subsection{Consequence 2}

The whole argument also works for $\mathbb{Q}$, or any field
$\mathbf{k}$ with $\mathbb{Q}\subset \mathbf{k}$, coefficients, and
we have:

\begin{theorem}\label {q}
Let $A=\mathbf{k}[x]/(x^{n+1})$, where $|x|=2m$, then as a
BV-algebra:
$$HH^*(A;A)=\mathbf{k}[x,u,t]/(x^{n+1},u^2,x^nt,ux^n),$$
where $|x|=-2m,\ |u|=-1\ and\ |t|=2mn+2(m-1)$, and
\begin{equation*}
\begin{split}
&\Delta(t^kx^l)=0\\
&\Delta(t^kux^l)=(-(k+1)n-k+l)t^kx^l.\\
\end{split}
\end{equation*}
\end{theorem}

In \cite{felix--thomas}, Y.Felix and J.Thomas proved that the BV-structures on $\mathbb{H}_*(LM;\mathbf{k})$
and on $HH^*(C^*(M;\mathbf{k});C^*(M;\mathbf{k}))$ can be naturally
identified. Precisely, they showed:

\textbf{Theorem}\cite{felix--thomas},\textbf{Theorem 1}: If $M$ is
$1$-connected and the field of coefficients has characteristic zero
then there exists a BV-structure on $HH^*(C^*(M);C^*(M))$ and an
isomorphism of BV-algebras $\mathbb{H}_*(LM)\cong
HH^*(C^*(M);C^*(M))$.

Therefore, theorem \ref {q} eventually calculates the BV-structure
of the rational loop homology,
$\mathbb{H}_*(L\mathbb{C}P^n;\mathbf{k})$ and
$\mathbb{H}_*(L\mathbb{H}P^n;\mathbf{k})$, of projective spaces.
We have:

\begin{theorem}\label{lh} Let $\mathbf{k}$ be any field containing $\mathbb{Q}$, then as BV-algebras:
$$\mathbb{H}_*(L\mathbb{C}P^n;\mathbf{k})=\mathbf{k}[x,u,t]/(x^{n+1},u^2,x^nt,ux^n),$$
with $|x|=-2,\ |u|=-1\ and\ |t|=2n$; and
$$\mathbb{H}_*(L\mathbb{H}P^n;\mathbf{k})=\mathbf{k}[x,u,t]/(x^{n+1},u^2,x^nt,ux^n),$$
with $|x|=-4,\ |u|=-1\ and\ |t|=4n+2$; and in both cases,

\begin{equation*}
\begin{split}
&\Delta(t^kx^l)=0\\
&\Delta(t^kux^l)=(-(k+1)n-k+l)t^kx^l.\\
\end{split}
\end{equation*}
\end{theorem}

\textit{Proof:} $\mathbf{K}P^n$'s are formal, i.e.,
$C^*(\mathbf{K}P^n)$ is quasi-isomorphic to $H^*(\mathbf{K}P^n)$.
Therefore, by the naturality of the Hochschild cohomology with
respect to quasi-isomorphisms \cite{felix--thomas--vigue-poirrier3},
\textbf{3} and \cite{felix--thomas},\textbf{Theorem 1}:
$$\mathbb{H}_*(L\mathbf{K}P^n;\mathbf{k})\cong HH^*(C^*(\mathbf{K}P^n;\mathbf{k});C^*(\mathbf{K}P^n;\mathbf{k}))\cong
HH^*(H^*(\mathbf{K}P^n;\mathbf{k});H^*(\mathbf{K}P^n;\mathbf{k}))$$
as BV-algebras. Theorem \ref {q} completes the proof, when $m=1$ and
$2$.\hfill$\square$

As another consequence of the main theorem, we have:
\begin{theorem}\label{tri}
Let $A=\mathbf{k}[x]/(x^{n+1})$, then the Lie bracket on the
negative cyclic cohomology $HC_-^*(A)$ is trivial.
\end{theorem}

It is a consequence of the following calculations and propositions.

Since $\varphi^*:C^*(A;A)\rightarrow P^*(A)$ is a chain map, it
takes the $\Delta$-operator of $C^*(A;A)$ onto the periodic
resolution $P^*(A)$ making it into a mixed complex $(P^*(A),d,B)$,
Moreover, $(P^*(A),d,B)$ is quasi-isomorphic to
$(C^*(A;A),\beta,\Delta)$ as mixed complexes, since $\varphi^*$ is a
quasi-isomorphism. Now we can use this simpler mixed complex to
calculate the negative cyclic cohomology of $A$. Precisely, as a
mixed complex, $B^{2q}=0$, $B^{2q+1}=B(q)\oplus 0$, where
$$B(q):\overline{A}\rightarrow\overline{A},\ \
B(q)(x^k)=(-q(n+1)-n+k)x^k.$$ Therefore, by tracing the defining
double complex $\textbf{BC}_-^{**}(A)$ of the negative cyclic
cohomology of $A$ given by the mixed complex $(P^*(A),d,B)$, we
have:



\begin{proposition} Let $\mathbf{k}$ be any field of characteristic
0, then
\begin{equation*}
\begin{split}
&HC_-^m(\mathbf{k}[x]/(x^{n+1}))=\mathbf{k},\ \ m=2q\ \
q=0,1,2,...\\
&HC_-^m(\mathbf{k}[x]/(x^{n+1}))=0,\ \ m=2q+1\ \ q=0,1,2,...\\
&HC_-^m(\mathbf{k}[x]/(x^{n+1}))=0,\ \ m=-2q\ \ q=1,2,...\\
&HC_-^m(\mathbf{k}[x]/(x^{n+1}))=\mathbf{k}^n,\ \ m=-2q+1\ \
q=1,2,...\\
\end{split}
\end{equation*}
\end{proposition}

Since in $\mathbf{k}$, all elements are invertible hence $B(q)$'s
are surjective.

Now it is time to compute the Lie bracket on $HC_-^*(A)$ defined in
\cite{menichi1}. First consider the Connes' long exact sequence for
negative cyclic cohomology
$$\cdots\rightarrow HH^n(A)\stackrel{I}{\rightarrow}HC_-^n(A)\rightarrow HC_-^{n+2}(A)\stackrel{\partial}{\rightarrow}HH^{n+1}\rightarrow\cdots.$$

By checking the definition of the connecting morphism carefully, we
have
$$\partial:HC_-^{n+2}(A)\rightarrow HH^{n+1}(A);\ \ \partial(a_{n+2},a_{n+4},...)=B(a_{n+2}).$$

For $a_i\in HC_-^{m_i}\ i=1,2$, $[a_1,a_2]$ is defined to be
$I(\partial(a_1)\cup\partial(a_2))$. Therefore, if $a\in
HC_-^{2q}(A)$, then $\partial (a)=B(a_{2q+2})=0$, hence
$[a,\_]=[\_,a]=0$. If $a_i\in HC_-^{2q_i+1}(A)$, then $[a_1,a_2]\in
HC_-^{2(q_1+q_2+1)}(A)=0$, hence trivial. Therefore, the Lie bracket
is trivial. \hfill$\square$


\subsection{Consequence 3}

Also, we have the similar result for $\mathbb{Z}_p$, and any field
$\mathbb{F}_p$ of characteristic $p$, coefficients, where $p$ is an
arbitrary prime number other than $2$.

When $n\neq kp-1$, all the argument in the $\mathbb{Z}$
coefficients case still works; and we have:

\begin{theorem}\label{p2}
Let $A=\mathbb{F}_p[x]/(x^{n+1})$, where $|x|=2m$ and $n\neq kp-1$,
then as a BV-algebra:
$$HH^*(A;A)=\mathbb{F}_p[x,u,t]/(x^{n+1},u^2,x^nt,ux^n),$$
where $|x|=-2m,\ |u|=-1\ and\ |t|=2mn+2(m-1)$, and
\begin{equation*}
\begin{split}
&\Delta(t^kx^l)=0\\
&\Delta(t^kux^l)=(-(k+1)n-k+l)t^kx^l.\\
\end{split}
\end{equation*}
\end{theorem}

When $n=kp-1$, we have:

\begin{theorem}\label{p}
Let $A=\mathbb{F}_p[x]/(x^{kp})$, where $|x|=2m$, as a BV-algebra:
$$HH^*(A;A)=\mathbb{F}_p[x,v,t]/(x^{kp},v^2),$$ where $|x|=-2m,\ |v|=2m-1\
and\ |t|=2mn+2(m-1)$, and
\begin{equation*}
\begin{split}
&\Delta(t^kx^l)=0\\
&\Delta(t^kvx^l)=lt^kx^{l-1}.\\
\end{split}
\end{equation*}
\end{theorem}

\textit{Proof:} We still use the $2$-periodical resolution and the
bridge $\varphi^*$ between $C^*(A;A)$ and $P^*(A)$. In this case,
$P^*(A)$ turns out to be

$$0\rightarrow A \stackrel{0}{\rightarrow} A
\stackrel{0}{\rightarrow} A \stackrel{0}{\rightarrow} A
\stackrel{0}{\rightarrow} A \stackrel{0}{\rightarrow} \cdots.$$

Let's take a close look at the following elements in this case.
\begin{equation*}
\begin{split}
&x\in P^0(A)=A,\ x(1)=x,\\
&v\in P^1(A)=A,\ v(1)=1,\ \text{and}\\
&t\in P^2(A)=A,\ t(1)=1.\\
\end{split}
\end{equation*}

We have $|x|=-2m,|v|=2m-1$,and $|t|=2mn+2(m-1)$; and these elements
correspond via $\varphi^*$ to
\begin{equation*}
\begin{split}
&\bar{x}\in Hom(\mathbb{F}_p,A),\ \bar{x}(1)=x,\\
&\bar{v}\in Hom(s\bar{A},A),\ \bar{v}(x^i)=ix^{i-1},\ \text{and}\\
&\bar{t}\in Hom(s\bar{A}^{\otimes 2},A),\ \bar{t}(x^i,
x^j)=x^{i+j-(n+1)}.\\
\end{split}
\end{equation*}

$\bar{v}$ is a cocycle because, if $i+j<n$,
$$\beta(\bar{v})(x^i,x^j)=x^i\bar{v}(x^j)-\bar{v}(x^{i+j})+\bar{v}(x^i)x^j=ix^{i+j-1}-(i+j)x^{i+j-1}+jx^{i+j-1}=0;$$
if $i+j=n+1$,
$$\beta(\bar{v})(x^i,x^j)=x^i\bar{v}(x^j)+\bar{v}(x^i)x^j=ix^{i+j-1}+jx^{i+j-1}=(n+1)x^n=kpx^n=0\in \mathbb{F}_p[x]/(x^{kp});$$
and if $i+j>n+1$, every term vanishes, so $\beta(\bar{v})=0$.

These elements represent the generators of $HH^*(A;A)$ as an
associative algebra. $v^2=0$ because $v^2(1\otimes
1)=\bar{v}^2(\varphi_2(1\otimes
1))=-\frac{n(n+1)}{2}x^{n-1}=-\frac{(kp-1)kp}{2}x^{n-1}$; if $k=2l$,
then $-\frac{(kp-1)kp}{2}=-(kp-1)lp=0\in \mathbb{F}_p$, and if
$k=2l+1$, then $kp-1=2l'$ since $p$ is odd, so
$-\frac{(kp-1)kp}{2}=l'kp=0 \in \mathbb{F}_p$. Then by checking the
formula of $\Delta$, we get the result.\hfill$\square$

In \cite {menichi2}, Luc Menichi conjectured that for any prime $p$,
the free loop space modulo $p$ of the complex projective space
$\mathbb{H}_*(L\mathbb{C}P^{p-1};\mathbb{Z}_p)$ is not isomorphic as
BV-algebras to the Hochschild cohomology
$HH^*(H^*(\mathbb{C}P^{p-1};\mathbb{Z}_p);H^*(\mathbb{C}P^{p-1};\mathbb{Z}_p))$.
He also pointed that: in \cite {Bokstedt--Ottosen},
M.B$\ddot{o}$kstedt and I.Ottosen have announced the computation of
BV-structure of $\mathbb{H}_*(L\mathbb{C}P^n;\mathbb{Z}_p)$.
Therefore, combining with theorem \ref {p}, this will give a
complete answer of Menichi's conjecture.

There is no essential difference between $p=2$ and $p=$ other prime
numbers, so we have:

\begin{theorem}

Let $A=\mathbb{F}_2[x]/(x^{n+1})$, where $|x|=2m$. If $n$ is even,
as a BV-algebra:

$$HH^*(A;A)=\mathbb{F}_2[x,u,t]/(x^{n+1},u^2,x^nt,ux^n),$$
where $|x|=-2m,\ |u|=-1\ and\ |t|=2mn+2(m-1)$, and
\begin{equation*}
\begin{split}
&\Delta(t^kx^l)=0\\
&\Delta(t^kux^l)=(-k+l)t^kx^l;\\
\end{split}
\end{equation*}
if $n$ is odd, as a BV-algebra:

$$HH^*(A;A)=\mathbb{F}_2[x,v,t]/(x^{n+1},v^2-\frac{n+1}{2}tx^{n-1}),$$ where $|x|=-2m,\ |v|=2m-1\
and\ |t|=2mn+2(m-1)$, and
\begin{equation*}
\begin{split}
&\Delta(t^kx^l)=0\\
&\Delta(t^kvx^l)=lt^kx^{l-1},\\
\end{split}
\end{equation*}
especially when $n=1$, as a BV-algebra:

$$HH^*(A;A)=\mathbb{F}_2[x,v,t]/(x^2,v^2-t)\cong
\Lambda[x]\otimes \mathbb{F}_2[v],$$ where $|x|=-2m$ and $|v|=2m-1$
and
\begin{equation*}
\begin{split}
&\Delta(v^k)=0\\
&\Delta(v^kx)=kv^{k-1}.\\
\end{split}
\end{equation*}

\end{theorem}

\textit{Proof:} Every step is exactly the same as the $\mathbb{F}_p$
coefficients case, except that when $n$ is odd, $v^2(1\otimes
1)=\frac{n(n+1)}{2}x^{n-1}=\frac{n+1}{2}x^{n-1}=\frac{n+1}{2}tx^{n-1}(1\otimes1).$\hfill$\square$

This recalculates the results in \cite {menichi2}, \cite
{westerland} and \cite {westerland2}.

\section*{Acknowledgement}

The author would like to thank Don Stanley for directing his
attention to this problem. It was also Don who suggests him to
construct the chain map $\varphi$ in the key step. He would also like to thank the referee for the careful revision and providing several instructive suggestions on improving this work, and J. Stasheff for showing interest to this work and warm encouragement to the author.

\bigskip
\noindent
Tian Yang\\
Department of Mathematics, Rutgers University\\
New Brunswick, NJ 08854, USA\\
(tianyang@math.rutgers.edu)

\end{document}